\documentclass[11pt]{article}
\usepackage{epsfig}

\def\be{\begin{equation}}
\def\ee{\end{equation}}
\def\bea{\begin{eqnarray}}
\def\eea{\end{eqnarray}}
\def\bes{\begin{eqnarray*}}
\def\ees{\end{eqnarray*}}

\def\nn{\nonumber}
\def\<{\langle}
\def\>{\rangle}
\def\lb{\label}
\def\bs{\setminus}

\def\R{{\bf R}}
\def\C{{\bf C}}
\def\Z{{\bf Z}}
\def\K{{\bf K}}
\def\N{{\bf N}}
\def\U{{\bf U}}

\def\Q{{\bf Q}}

\def\aa{{\alpha}}

\def\ga{{\gamma}}

\def\ka{{\kappa}}
\def\th{{\theta}}

\def\Om{{\Omega}}
\def\ep{{\epsilon}}
\def\lm{{\lambda}}
\def\Lm{{\Lambda}}

\def\sg{{\sigma}}
\def\dm{{\diamond}}

\def\K{{\cal K}}

\def\Sp{{\rm Sp}}

\def\ol{\overline}

\def\ol#1{\overline{#1}}  %overline in math mode

\def\hb{\vrule height0.18cm width0.14cm $\,$}

\def\ol#1{\overline{#1}}  %overline in math mode

\title{Multiple closed geodesics on bumpy Finsler $n$-spheres}

\author{Huagui Duan$^{1}$,\thanks{Partially supported by NNSF and RFDP
of MOE of China. E-mail: duanhuagui@163.com}
\qquad Yiming Long$^{1,2}$\thanks{Partially supported by the 973 Program of
MOST, Yangzi River Professorship, NNSF, MCME, RFDP, LPMC of MOE of
China, S. S. Chern Foundation, and Nankai University.
E-mail: longym@nankai.edu.cn }\\ \\
$^{1}$ Chern Institute of Mathematics\\
$^{2}$ Key Lab of Pure Mathematics and Combinatorics of Ministry of
Education\\ Nankai University, Tianjin 300071\\ The People's
Republic of China\\ }
\date{}

\begin{document}

\maketitle

\begin{abstract}
{\it  In this paper we prove that for every bumpy Finsler metric
$F$ on every rationally homological $n$-dimensional sphere $S^n$
with $n\ge 2$, there exist always at least two distinct prime
closed geodesics.}
\end{abstract}

\renewcommand{\theequation}{\thesection.\arabic{equation}}
\renewcommand{\thefigure}{\thesection.\arabic{figure}}

\setcounter{equation}{0}%\setcounter{figure}{0}
\section{Introduction and the main result}%{Section 1}

Let us recall firstly the definition of the Finsler metric.

{\bf Definition 1.1} (cf. \cite{BCS1} and \cite{She1}) {\it Let $M$
be a finite dimensional manifold and $TM$ be its tangent bundle. A
function $F:TM\rightarrow [0,+\infty)$ is a Finsler metric if it
satisfies the following properties:

$(F_1)$ $F$ is $C^{\infty}$ on $TM\backslash\{0\}$.

$(F_2)$ $F(\lm y)=\lm F(y)$ for all $\lm>0$ and $y\in TM$.

$(F_3)$ For any $y\in TM\backslash \{0\}$, the symmetric bilinear
form $g_y$ on $TM$ is positive definite, where
$$ g_y(u,v)=\frac{1}{2}
  \frac{\partial^2}{\partial s \partial t}[F^2(y+su+tv)]|_{s=t=0}. $$
The pair $(M,F)$ is called a Finsler manifold. A Finsler metric $F$ is
reversible if $F(-v)=F(v)$ for all $v\in TM$.}

For the definition of closed geodesics on a Finsler manifoid, we
refer readers to \cite{BCS1} and \cite{She1}. As usual, on any
Finsler manifold $M=(M,F)$ a closed geodesic $c:S^1=\R/\Z\to M$
is {\it prime}, if it is not a multiple covering (i.e., iteration)
of any other closed geodesics. Here the $m$-th iteration $c^m$ of
$c$ is defined by $c^m(t)=c(mt)$ for $m\in\N$. The inverse curve
$c^{-1}$ of $c$ is defined by $c^{-1}(t)=c(1-t)$ for $t\in \R$. We
call two prime closed geodesics $c$ and $d$ {\it distinct} if there
is no $\th\in (0,1)$ such that $c(t)=d(t+\th)$. We shall omit the
word "distinct" for short when we talk about more than one prime
closed geodesics. A closed geodesic $c$ on $(M,F)$ is non-degenerate,
if its linearized Poincar\'e map $P_c$ has no eigenvalue $1$. A
Finsler metric $F$ on $M$ is {\it bumpy} if all closed geodesics and
their iterates on $(M,F)$ are non-degenerate.

In recent years, geodesics and closed geodesics on Finsler manifolds
have got more attentions. We refer readers to \cite{BRS1} of D. Bao,
C. Robles and Z. Shen, \cite{Rob1} of C. Robles, and \cite{Lon3}
of Y. Long and the references therein for recent progress in this area.

Note that by the classical theorem of Lyusternik-Fet \cite{LyF1} in 1951,
there exists at least one closed geodesic on every compact Riemannian
manifold. Because the proof is variational, this result works also for
compact Finsler manifolds. In \cite{Rad3} of 2005, H.-B. Rademacher
obtained existence of closed geodesics on $n$-dimensional Finsler spheres
under pinching conditions which generalizes results in \cite{BTZ1} and
\cite{BTZ2} of W. Ballmann, W. Thorbergsson and W. Ziller in 1982-83
on Riemannian manifolds.

We are only aware of a few results on the existence of multiple closed
geodesics on Finsler spheres without pinching conditions. In \cite{Fet1}
of 1965, A. I. Fet proved that there exist at least two distinct closed
geodesics on every reversible bumpy Finsler manifold $(M,F)$. In \cite{Rad1}
of 1989, H.-B. Rademacher proved that there exist at least two elliptic closed
geodesics on every bumpy Finsler 2-sphere. In \cite{HWZ1} of 2003, H. Hofer,
K. Wysocki and E. Zehnder proved that there exist either two or
infinitely many distinct closed geodesics on every bumpy Finsler $2$-sphere
if the stable and unstable manifolds of every hyperbolic closed geodesics
intersect transversally. In \cite{BaL1} of 2005, V. Bangert and Y. Long
proved that there exist at least two distinct prime closed geodesics on
every Finsler $2$-sphere $(S^2,F)$.

The aim of this paper is to prove the following main result,
specially for bumpy irreversible Finsler rationally homological
$n$-spheres without pinching conditions.

{\bf Theorem 1.2.} {\it For every bumpy Finsler metric $F$ on every rationally
homological $n$-sphere $S^n$ with $n\ge 2$, there exist at least two distinct
prime closed geodesics. }

Note that our proof of Theorem 1.2 uses only the $\Q$-homological properties
of the Finsler manifold, thus we shall carry out our proof of this theorem
below just for $n$-dimensional spheres.

In this paper, let $\N$, $\N_0$, $\Z$, $\Q$, $\R$, and $\C$ denote
the sets of positive integers, non-negative integers, rational
numbers, real numbers and complex numbers respectively. We denote
by $[a]=\max\{k\in\Z\,|\,k\le a\}$ for any $a\in\R$. We use only
singular homology modules with $\Q$-coefficients.

\setcounter{equation}{0}
\section{Critical modules of iterations of closed geodesics}%{Section 2}

Let $M=(M,F)$ be a compact Finsler manifold $(M,F)$, the space
$\Lambda=\Lambda M$ of $H^1$-maps $\gamma:S^1\rightarrow M$ has a
natural structure of Riemannian Hilbert manifolds on which the
group $S^1=\R/\Z$ acts continuously by isometries, cf.
\cite{Kli2}, Chapters 1 and 2. This action is defined by
$(s\cdot\gamma)(t)=\gamma(t+s)$ for all $\gamma\in\Lm$ and $s,
t\in S^1$. For any $\gamma\in\Lambda$, the energy functional is
defined by
\be E(\gamma)=\frac{1}{2}\int_{S^1}F(\gamma(t),\dot{\gamma}(t))^2dt.
\lb{2.1}\ee
It is of class $C^{1,1}$ (cf. \cite{Mer1}) and
invariant under the $S^1$-action. The critical points of $E$ of
positive energies are precisely the closed geodesics
$\gamma:S^1\to M$. The index form of the functional $E$ is well
defined along any closed geodesic $c$ on $M$, which we denote by
$E''(c)$ (cf. \cite{She1}). As usual, we denote by $i(c)$ and
$\nu(c)$ the Morse index and nullity of $E$ at $c$. In the
following, we denote by
\be \Lm^\kappa=\{d\in \Lm\;|\;E(d)\le\kappa\},\quad \Lm^{\kappa-}=\{d\in \Lm\;|\; E(d)<\kappa\},
  \quad \forall \kappa\ge 0. \lb{2.2}\ee

For $m\in\N$ we denote the $m$-fold iteration map
$\phi_m:\Lambda\rightarrow\Lambda$ by $\phi_m(\ga)(t)=\ga(mt)$,
for all $\,\ga\in\Lm, t\in S^1$, as well as $\ga^m=\phi_m(\gamma)$.
For a closed geodesic $c$, recall that the mean index $\hat{i}(c)$ is
defined by
\be \hat{i}(c)=\lim_{m\rightarrow\infty}\frac{i(c^m)}{m}. \lb{2.3}\ee

If $\gamma\in\Lambda$ is not constant then the multiplicity $m(\gamma)$
of $\gamma$ is the order of the isotropy group
$\{s\in S^1\mid s\cdot\gamma=\gamma\}$. If $m(\gamma)=1$ then $\gamma$ is
prime. Hence $m(\gamma)=m$ if and only if there exists a prime curve
$\tilde{\gamma}\in\Lambda$ such that $\gamma=\tilde{\gamma}^m$.

For a closed geodesic $c$ we set $ \Lm(c)=\{\ga\in\Lm\mid
E(\ga)<E(c)\}. $ If $A\subseteq\Lm$ is invariant under the action of some
subgroup $\Gamma$ of $S^1$, we denote by $A/\Gamma$ the quotient
space of $A$ module the action of $\Gamma$.

Using singular homology with rational coefficients we will
consider the following critical $\Q$-module of a closed geodesic
$c\in\Lambda$: \be \overline{C}_*(E,c)
   = H_*\left((\Lm(c)\cup S^1\cdot c)/S^1,\Lm(c)/S^1\right). \lb{2.4}\ee

In order to apply the results of D. Gromoll and W. Meyer in \cite{GrM1} and
\cite{GrM2}, following \cite{Rad2}, Section 6.2, we introduce finite-dimensional
approximations to $\Lambda$. We choose an arbitrary energy value $a>0$ and
$k\in\N$ such that every  geodesic segment of length $<\sqrt{2a/k}$ is minimal.
Then
$$ \Lm(k,a)=\left\{\ga\in\Lm \mid E(\ga)<a \mbox{ and }
    \ga|_{[i/k,(i+1)/k]}\mbox{ is a geodesic segment for }i=0,\ldots,k-1\right\} $$
is a $(k\cdot\dim M)$-dimensional submanifold of $\Lambda$
consisting of closed geodesic polygons with $k$ vertices. The set
$\Lambda(k,a)$ is invariant under the action of the subgroup $\Z_k$ of $S^1$.
Closed geodesics in $\Lambda^{a-}=\{\gamma\in\Lambda\mid
E(\gamma)<a\}$ are precisely the critical points of
$E|_{\Lm(k,a)}$, and for every closed geodesic $c\in\Lm(k,a)$ the
index of $(E|_{\Lm(k,a)})''(c)$ equals $i(c)$ and the null space
of $(E|_{\Lm(k,a)})''(c)$ coincides with the nullspace of
$E''(c)$, cf. \cite{Rad2}, p.51.

We call a closed geodesic satisfying the isolation condition, if
the following holds:

{\bf (Iso) The orbit $S^1\cdot c^m$ is an isolated critical orbit
of $E$ for all $m\in\N$. }

Since our aim is to prove the existence of more than one closed
geodesic for every bumpy Finsler metric on $S^n$, the condition
{\bf (Iso)} does not restrict generality.

Now we can apply the results by D. Gromoll and W. Meyer
\cite{GrM1} to a given closed geodesic $c$ satisfying (Iso). If
$m=m(c)$ is the multiplicity of $c$, we choose a
finite-dimensional approximation $\Lm(k,a)\subseteq\Lm$ containing
$c$ such that $m$ divides $k$. Then the isotropy subgroup
$\Z_m\subseteq S^1$ of $c$ acts on $\Lm(k,a)$ by isometries. Recall
that the $\Z_m$-action is defined by $\frac{i}{m}\cdot g(t)=g(t+\frac{i}{m})$
for all $g\in\Lm(k,a)$ and $\frac{i}{m}\in\Z_m$ with $1\le i\le m$.
Let $D$ be a $\Z_m$-invariant local hypersurface transverse to
$S^1\cdot c$ in $c\in \Lm(k,a)$. Such a $D$ can be obtained by
applying the exponential map of $\Lm(k,a)$ at $c$ to the normal
space to $S^1\cdot c$ at $c$. We denote by
\be T_cD=V_+\oplus V_-\oplus V_0,   \lb{2.5}\ee
the orthogonal decomposition of $T_cD$ into the positive, negative and null
eigenspace of the endomorphism of $T_cD$ associated to $(E|_D)''(c)$ by the
Riemannian metric. In particular, we have $\dim V_-=i(c)$ and
$\dim V_0=\nu(c)$. According to \cite{GrM1}, Lemma 1, for every such a $D$
there exist balls $B_+\subseteq V_+, B_-\subseteq V_-$ and $B_0\subseteq V_0$
centered at the origins, a diffeomorphism
$$  \psi:B=B_+\times B_-\times B_0\rightarrow\psi(B_+\times B_-\times B_0)\subseteq D  $$
with $\psi(0)=c$, $\psi_{*0}$ preserving the splitting
(\ref{2.5}), and a smooth function $f:B_0\rightarrow\R$ satisfying
\bea
&& f^\prime(0)=0 \quad \mbox{and} \quad f^{\prime\prime}(0)=0,  \lb{2.6}\\
&& E\circ\psi(x_+,x_-,x_0)=|x_+|^2 - |x_-|^2 + f(x_0),
\lb{2.7}\eea
for $(x_+, x_-, x_0) \in B_+\times B_-\times B_0$. Since the $\Z_m$-action is
isometric and $E$ is $\Z_m$-invariant, the tangential map
$(\frac{i}{m}|_D)_{*c}$ of $\frac{i}{m}\in\Z_m$ restricted to $D$ at $c$ preserves the
above splitting (\ref{2.5}). It follows from the construction of $\psi$ that
$\psi$ is equivariant with respect to the $\Z_m$-action, i.e.,
$\frac{i}{m}\cdot\psi=\psi\circ(\frac{i}{m}|_D)_{*c}\cdot\,$ for $\frac{i}{m}\in\Z_m$,
cf. p.501 of \cite{GrM2}.

As in \cite{GrM1} and \cite{GrM2}, we call
$N=\{\psi(0,0,x_0)\,|\,x_0\in B_0\}$ a local characteristic
manifold at $c$, $U=\{\psi(0,x_-,0)\,|\,x_-\in B_-\}$ a local
negative disk at $c$. Note that $N$ and $U$ are $\Z_m$-invariant.
It follows from (\ref{2.7}) that $c$ is an isolated critical point
of $E|_N$. We set $N^-=N\cap\Lm(c)$,
$U^-=U\cap\Lm(c)=U\setminus\{c\}$ and $D^-=D\cap\Lm(c)$. Using
(\ref{2.7}), the fact that $c$ is an isolated critical point of
$E|_N$, and the K\"unneth formula, we obtain \bea &&
H_*(D^-\cup\{c\},D^-)=
      H_*(U^-\cup\{c\},U^-) \otimes H_*(N^-\cup\{c\},N^-),
      \lb{2.8}\\[5pt]
&& H_q(U^-\cup\{c\}, U^-) = H_q(U, U\setminus\{c\})
    = \left\{\matrix{\Q, & {\rm if\;}q=i(c), \cr
                      0, & {\rm otherwise}, \cr}\right.,  \lb{2.9}\eea
cf. \cite{Rad2}, Lemma 6.4 and its proof. As in p.59 of
\cite{Rad2}, for all $m\in\N$, let respectively
\be H_{\ast}(X,A)^{\pm\Z_m}
   = \{[\xi]\in H_{\ast}(X,A)\,|\,T_{\ast}[\xi]=\pm [\xi]\}, \lb{2.10}\ee
where $T$ is a generator of the $\Z_m$ action.

Now we have the following Propositions.

{\bf Proposition 2.1.} (cf. Satz 6.11 of \cite{Rad2} ) {\it Let
$c$ be a prime closed geodesic on a Finsler manifold $(M,F)$
satisfying (Iso). Then we have \bea \overline{C}_q( E,c^m)
&\equiv& H_q\left((\Lm(c^m)\cup S^1\cdot c^m)/S^1, \Lm(c^m)/S^1\right)\nn\\
&=& \left(H_{i(c^m)}(U_{c^m}^-\cup\{c^m\},U_{c^m}^-)
    \otimes H_{q-i(c^m)}(N_{c^m}^-\cup\{c^m\},N_{c^m}^-)\right)^{+\Z_m} \nn
\eea

(i) When $\nu(c^m)=0$, there holds
$$ \overline{C}_q( E,c^m) = \left\{\matrix{
     \Q, &\quad {\it if}\;\; i(c^m)-i(c)\in 2\Z\;\;{\it and}\;\;
                   q=i(c^m),\;  \cr
     0, &\quad {\it otherwise}. \cr}\right.  $$

(ii) When $\nu(c^m)>0$, there holds
$$ \overline{C}_q( E,c^m) =
    H_{q-i(c^m)}(N_{c^m}^-\cup\{c^m\}, N_{c^m}^-)^{\ep(c^m)\Z_m}, $$
where $\ep(c^m)=(-1)^{i(c^m)-i(c)}$.}

We need the following

{\bf Definition 2.2.} (cf. \cite{Rad2}, \cite{BaL1}, \cite{LoW1}) {\it Suppose
$c$ is a closed geodesic of multiplicity $m(c)=m$ satisfying (Iso). If $N$ is
a local characteristic manifold at $c$, $N^- = N\cap\Lm(c)$ and $j\in\Z$, we define}
\bea
k_j(c) &\equiv& \dim\, H_j( N^-\cup\{c\},N^-),    \nn\\
k_j^{\pm 1}(c) &\equiv& \dim\, H_j(N^-\cup\{c\},N^- )^{\pm\Z_m}. \nn \eea

Clearly the integers $k_j(c)$ and $k_j^{\pm 1}(c)$ equal to $0$ when $j<0$ or
$j>\nu(c)$, and can take only values $0$ or $1$ when $j=0$ or $j=\nu(c)$.

{\bf Proposition 2.3.} (cf. Satz 6.13 of \cite{Rad2}, \cite{BaL1}, \cite{LoW1})
{\it Let $c$ be a prime closed geodesic satisfying (Iso).

(i) There holds $0\le k_j^{\pm 1}(c^m) \le k_j(c^m)$ for all $m\in\N$ and $j\in\Z$.

(ii) For any $m\in\N$, there hold $k_0^{+1}(c^m) = k_0(c^m)$ and $k_0^{-1}(c^m) = 0$.

(iii) In particular, if $c^m$ is non-degenerate, i.e. $\nu(c^m)=0$, then
$k_0^{+1}(c^m) = k_0(c^m)=1$ and $k_0^{-1}(c^m) = 0$ hold. }

\setcounter{equation}{0}
\section{The structure of $H_{\ast}(\ol{\Lm} S^n,\ol{\Lm}^0S^n;\Q)$ }%{Section 3}

In this section, we briefly describe the relative homological
structure of the quotient space $\overline{\Lm}\equiv\overline{\Lm} S^n=\Lm S^n/S^1$.
Here we have
$\ol{\Lm}^0=\ol{\Lambda}^0S^n =\{{\rm constant\;point\;curves\;in\;}S^n\}\cong S^n$.

Let $(X,Y)$ be a space pair such that the Betti numbers
$b_i=b_i(X,Y)=\dim H_i(X,Y;\Q)$ are finite for all $i\in \Z$. As
usual the {\it Poincar\'e series} of $(X,Y)$ is defined by the
formal power series $P(X, Y)=\sum_{i=0}^{\infty}b_it^i$. We need the
following well known results on Betti numbers and the Morse inequality.

{\bf Theorem 3.1.} (cf. Theorem 2.4 and Remark 2.5 of \cite{Rad1}) {\it

(i) When $n\in 2\N$, we have
\bea P(\ol{\Lm}S^n,\ol{\Lm}^0S^n)(t)
&=& t^{n-1}\left(\frac{1}{1-t^2}+\frac{t^{2n-2}}{1-t^{2n-2}}\right) \nn\\
&=& (t^{(n-1)}+t^{(n+1)}+t^{(n+3)}+\cdots)+(t^{3(n-1)}+t^{5(n-1)}+t^{7(n-1)}+\cdots), \nn\eea
which yields
\be b_q  \equiv \dim H_q(\ol{\Lm} S^n,\ol{\Lm}^0 S^n )
 = \;\;\left\{\matrix{
    2,&\quad if\quad q\in\K\equiv\{k(n-1)\,|\,3\le k\in(2\N+1)\},  \cr
   1,&\quad if\quad q\in\{(n-1)+2k\,|\,k\in\N_0\}\setminus\K,  \cr
    0, &\quad  otherwise. \cr}\right. \lb{3.1}\ee

(ii) When $n\in (2\N+1)$, we have
\bea P(\ol{\Lm}S^n,\ol{\Lm}^0S^n)(t)
&=& t^{n-1}\left(\frac{1}{1-t^2}+\frac{t^{n-1}}{1-t^{n-1}}\right) \nn\\
&=& (t^{(n-1)}+t^{(n+1)}+t^{(n+3)}+\cdots)+(t^{2(n-1)}+t^{3(n-1)}+t^{4(n-1)}+\cdots), \nn\eea
which yields}
\be b_q  \equiv \dim H_q(\ol{\Lm} S^n,\ol{\Lm}^0 S^n )
 = \;\;\left\{\matrix{
    2,&\quad if\quad q\in\K\equiv\{k(n-1)\,| \,2\le k\in\N\},  \cr
    1,&\quad if\quad q\in\{(n-1)+2k\,|\,k\in\N_0\}\setminus\K,  \cr
    0, &\quad  otherwise. \cr}\right. \lb{3.2}\ee

{\bf Theorem 3.2.} (cf. Theorem I.4.3 of \cite{Cha1}, Theorem 6.1
of \cite{Rad2}) {\it Suppose that there exist only finitely many
prime closed geodesics $\{c_j\}_{1\le j\le k}$ on a Finsler
$n$-sphere $(S^n, F)$. Set
$$ M_q =\sum_{1\le j\le k,\; m\ge 1}\dim{\ol{C}}_q(E, c^m_j), \quad \forall q\in\Z. $$
Then for every integer $q\ge 0$ there holds }
\bea
M_q - M_{q-1} + \cdots +(-1)^{q}M_0
&\ge& b_q - b_{q-1}+ \cdots + (-1)^{q}b_0, \lb{3.3}\\
M_q &\ge& b_q. \lb{3.4}\eea

\setcounter{equation}{0}
\section{Classification of closed geodesics on bumpy Finsler manifolds}%{Section 4}

Let $c$ be a closed geodesic on a Finsler manifold $(M,F)$. Denote the
linearized Poincar\'e map of $c$ by $P_c$. By \cite{LLo1} in 2002 of C. Liu and
Y. Long (cf. Chapter 12 of \cite{Lon2}), the index iteration formulae in \cite{Lon1}
work for Morse indices of iterated closed geodesics on Riemannian as well as Finsler
manifolds. We call a closed geodesic $c$ is {\it completely non-degenerate}, if
$c^m$ is non-degenerate for all $m\in\N$. When the Finsler metric $F$ is bumpy, every
closed geodesic $c$ on $(M,F)$ is completely non-degenerate. Thus by Theorems 8.1.4
to 8.1.7, 8.2.3 and 8.2.4, and 8.3.1 of \cite{Lon2}, in the basic normal form
decomposition of the symplectic matrix $P_c$ (cf. Theorem 1.8.10 of \cite{Lon2}) there
can exist only basic normal forms like $H(d)$ with $d\in\R\bs\{0,\pm 1\}$, $R(\th)$
and $N(\aa,B)$ with $\th/\pi$ and $\aa/\pi$ being irrational (cf. notation below).
Therefore according to the iteration formula of Morse indices, completely non-degenerate
closed geodesics on a Finsler manifold $(M,F)$ can be classified into the following 5
cases NCG-1 to NCG-5.

To introduce this classification, we need some notations from \cite{Lon2}. Given any
two real matrices of the square block form
$$ M_1=\left(\matrix{A_1 & B_1\cr C_1 & D_1\cr}\right)_{2i\times
2i},\qquad M_2=\left(\matrix{A_2 & B_2\cr C_2 &
D_2\cr}\right)_{2j\times 2j},$$ the $\diamond$-sum of $M_1$ and
$M_2$ is defined by the $2(i+j)\times2(i+j)$ matrix
$$ M_1\diamond M_2=\left(\matrix{A_1 & 0 & B_1 & 0 \cr
                                   0 & A_2 & 0& B_2\cr
                                   C_1 & 0 & D_1 & 0 \cr
                                   0 & C_2 & 0 & D_2}\right).$$
For convenience, we denote by
$N(\aa,B)^{\dm r}\equiv N({\aa}_1,{B}_1)\dm\cdots\dm N({\aa}_r,{B}_r)$,
where $\aa=(\aa_1,\ldots,\aa_r)$ and $B=(B_1,\ldots,B_r)$ for some
$0\le r\le [\frac{n-1}{2}]$. If $r=0$ in the following, it means
that no such a term $N(\aa,B)^{\dm r}$ appears. Here as in \cite{Lon2} we set
\bea
&& N({\aa}_i,{B_i})=\left(\matrix{R(\aa_i) & B_i      \cr
                                 0        & R(\aa_i) \cr}\right), \nn\\
&& R(\aa_i)=\left(\matrix{\cos\aa_i & -\sin\aa_i\cr
                              \sin\aa_i &  \cos\aa_i\cr}\right), \quad
B_i=\left(\matrix{b_{i1} & b_{i2}\cr b_{i3} & b_{i4}\cr}\right),
\nn\eea where $\aa_i/\pi\in (0, 2)\setminus(\Q\cup\{1\})$,
$(b_{i1},b_{i2},b_{i3},b_{i4})\in\R^4$ for $1\le i\le r$. We denote also
$H(d)=\left(\matrix{d & 0  \cr
                    0 & 1/d\cr}\right)$ with $d\in \R\bs\{0,\pm 1\}$.

The homotopy set $\Omega(M)$ of $M$ in the symplectic group $\Sp(2n)$ was
studied in \cite{Lon1} which is defined by
$$ \Om(M)=\{N\in \Sp(2n)\,|\,\sg(N)\cap\U=\sg(M)\cap\U\equiv\Gamma\;\mbox{and}
        \;\nu_{\omega}(N)=\nu_{\omega}(M)\, \forall\omega\in\Gamma\}, $$
where $\sigma(M)$ denotes the spectrum of $M$,
$\nu_{\omega}(M)\equiv\dim_{\C}\ker_{\C}(M-\omega I)$ for all
$\omega\in\U$, and $\U=\{z\in\C\;|\;|z|=1\}$. Let $\Omega^0(M)$ denote the path
connected component of $\Omega(M)$ containing $M$ (cf. p.38 of \cite{Lon2}).

By Theorems 8.2.3 to 8.2.4 of \cite{Lon2}, the Morse indices of iterates
of a completely non-degenerate closed geodesic $c$ with $P_c=N(\aa_i,B_i)$
satisfy the same formula
$$ i(c)=2p\ \mbox{for some}\ p\in\N_0\qquad\mbox{and}\qquad
i(c^m) = 2mp , \quad \nu(c^m)=0, \ \forall \ m\ge 1. $$
Hence by Theorems 8.1.4 to 8.1.7 and Theorem 8.3.1 of \cite{Lon2}, we have
the following classification of completely non-degenerate closed geodesics $c$
on a Finsler $n$-dimensional manifold, i.e., there exists a path
$f_c\in C([0,1],\Omega^0(P_c))$ such that $f_c(0)=P_c$ and $f_c(1)$ have the
following forms:

{\bf NCG-1.} {\it $f_c(1)= {N}(\aa,B)^{\diamond r}\diamond
R(\th_1)\diamond\cdots\diamond R(\th_{n-2r-1})$.}

In this case, by Theorem 8.3.1 of \cite{Lon2}, we have
$i(c)=2p+(n-2r-1)$ for some $p\in\Z$ such that $i(c)\ge0$, and \be
i(c^m) = 2mp
+2\sum_{i=1}^{n-2r-1}\left[\frac{m\th_i}{2\pi}\right]+(n-2r-1) ,
\quad \nu(c^m)=0, \qquad \forall \, m\ge 1. \lb{4.1}\ee

{\bf NCG-2.} {\it $f_c(1)=N(\aa,B)^{\diamond r}\diamond
R(\th_1)\diamond\cdots\diamond R(\th_{k})\diamond
H(d_{k+1})\diamond \cdots \diamond H(d_{n-2r-1})$ with $k\in 2\N$
and $2\le k\le n-2r-2$.}

In this case, by Theorem 8.3.1 of \cite{Lon2}, we have $i(c)=p$
for some $p\in\N_0$, and \be i(c^m) = m(p-k)
+2\sum_{i=1}^{k}\left[\frac{m\th_i}{2\pi}\right]+k , \quad
\nu(c^m)=0, \qquad \forall \ m\ge 1.\lb{4.2} \ee

{\bf NCG-3.} {\it $f_c(1)=N(\aa,B)^{\diamond r}\diamond
R(\th_1)\diamond\cdots\diamond R(\th_{k})\diamond
H(d_{k+1})\diamond \cdots \diamond H(d_{n-2r-1})$ with $k\in
(2\N-1)$ and $3\le k\le n-2r-2$.}

In this case, by Theorem 8.3.1 of \cite{Lon2}, we have $i(c)=p$
for some $p\in\N_0$, and \be i(c^m) =
m(p-k)+2\sum_{i=1}^{k}\left[\frac{m\th_i}{2\pi}\right]+k , \quad
\nu(c^m)=0, \qquad \forall \ m\ge 1.\lb{4.3} \ee

{\bf NCG-4.} {\it $f_c(1)=N(\aa,B)^{\diamond r}\diamond
R(\th_1)\diamond H(d_{2})\diamond \cdots \diamond
H(d_{n-2r-1})$.}

In this case, by Theorem 8.3.1 of \cite{Lon2}, we have $i(c)=p$
for some $p\in\N_0$, and \be i(c^m) =
m(p-1)+2\left[\frac{m\th_1}{2\pi}\right]+1 , \quad \nu(c^m)=0,
\qquad \forall \ m\ge 1.\lb{4.4} \ee

{\bf NCG-5.} {\it $f_c(1)=N(\aa,B)^{\diamond r}\diamond
H(d_1)\diamond \cdots \diamond H(d_{n-2r-1})$.}

In this case, by Theorem 8.3.1 of \cite{Lon2}, we have $i(c)=p$ for some
$p\in\N_0$, and
\be i(c^m) = mp, \quad \nu(c^m)=0, \qquad \forall \ m\ge 1. \lb{4.5}\ee

\setcounter{equation}{0}
\section{A mean index identity}%{Section 5}

We need a notation from \cite{LoW1}.

{\bf Definition 5.1.} {\it Let $c$ be a completely non-degenerate prime
closed geodesic on $(S^n,F)$. For each $m\in\N$, the critical type
numbers of $c^m$ is defined by
\be K(c^m) \equiv(k_0^\ep(c^m), k_1^\ep(c^m), \cdots,
   k_{n}^\ep(c^m)) =(k_0^{\ep }(c^m), 0, \cdots, 0), \lb{5.1}\ee
where $\ep=\ep(c^m)=(-1)^{i(c^m)-i(c)}$. Note that only
$k_0^\ep(c^m)$ may be non-zero for $m\ge 1$ by Definition 2.2 and
Proposition 2.3. We call a completely non-degenerate prime closed
geodesic $c$ homologically invisible if $k_0^\ep(c^m)=0$ for all
$m\in\N$, or homologically visible otherwise.}

{\bf Lemma 5.2.} {\it  Let $c$ be a completely non-degenerate prime closed
geodesic on a Finsler $n$-sphere $(S^n,F)$. Then there exist a minimal
integer $N\in\N$ such that $K(c^m)=K(c^{m+N})$ for all $m\in\N$.
According to the classification in Section 4, we have }
\bea
N&=&1, \qquad {\it if\;}c\;\mbox{\it belongs to NCG-1;}\nn\\
N&=&\left\{\matrix{
    1, \quad {\rm if}\;p\;{\rm is\;even},\cr
    2, \quad {\rm if}\;p\;{\rm is\;odd},\cr}\right.
    \quad {\it if\;}c\;\mbox{\it belongs to NCG-2 or NCG-5;}\nn\\
N&=&\left\{\matrix{
    2, \quad {\rm if}\;p\;{\rm is\;even},\cr
    1, \quad {\rm if}\;p\;{\rm is\;odd},\cr}\right.
    \quad {\it if\;}c\;\mbox{\it belongs to NCG-3 or NCG-4.} \nn
\eea

{\bf Proof.} In fact, $N$ depends only on the parity of
$i(c^m)-i(c)$ for any $m\in\N$ by Proposition 2.3. More precisely,
$$ N=\left\{\matrix{
    1, &\quad {\rm if}\;i(c^m)-i(c)\;{\rm is\;even\;for\;any}\;m\in\N,\cr
    2, &\quad {\rm otherwise}.\cr}\right. $$
By the classification of Section 4, we have the following details.
In NCG-1, $i(c^m)-i(c)$ is even. In NCG-2 and NCG-3,
$i(c^m)-i(c)=(m-1)(p-k)$ mod $2$. In NCG-4, $i(c^m)-i(c)=(m-1)(p-1)$ mod $2$.
In NCG-5, $i(c^m)-i(c)=(m-1)p$. Therefore Lemma 5.2 follows. \hfill\hb

Suppose that there exist only finitely many completely
non-degenerate prime closed geodesics $\{c_j\}_{1\le j\le k}$ for
$1\le j\le k$ on a bumpy Finsler $n$-sphere $(S^n, F)$. The
{\it Morse series} $M(t)$ of the energy functional $E$ on the space
$(\Lambda S^n/S^1, \Lambda^0 S^n/S^1)$ is defined by
$$ M(t)=\sum_{q\ge 0,\; m\neq 0 \atop 1\le j\le k}\dim\ol{C}_q(E, c^m_j)t^q. $$
Then it yields a formal power series $Q(t)=\sum_{i=0}^\infty
q_it^i$ with nonnegative integer coefficients $q_i$ such that \be
M(t)=P(\Lm S^n/S^1,\Lm^0 S^n/S^1)(t)+(1+t)Q(t). \lb{5.2}\ee For a
formal power series $R(t)=\sum_{i=0}^\infty r_it^i$, we denote by
$R^n(t)=\sum_{i=0}^n r_i t^i$ for $n\in\N$ the corresponding
truncated polynomials. Using this notation, (\ref{5.2}) becomes
\be (-1)^mq_m=M^m(-1)-P^m(-1) \qquad \forall m\in\N. \lb{5.3}\ee
By Satz 7.8 of \cite{Rad2} we have specially for spheres:
\be \lim_{m\to\infty}\frac{1}{m}P^m(\Lm S^n/S^1,\Lm^0 S^n/S^1)(-1) =
\left\{\matrix{-\frac{n}{2(n-1)}, & {\rm if\;}$n$\; {\rm is\; even}, \cr
              \frac{n+1}{2(n-1)}, & {\rm if\;}$n$\; {\rm is\; odd}. \cr}\right.
   \lb{5.4}\ee

A general version of the following mean index identity was proved in Theorem 3
in \cite{Rad1} and \cite{Rad2} of H.-B. Rademacher. Our following theorem gives
more precise coefficients in the identity than those in \cite{Rad1} and \cite{Rad2}.
This more precise information is crucial in the proof of our main Theorem 1.2
later.

{\bf Theorem 5.3.} {\it Suppose that there exist only finitely many homologically
visible prime closed geodesics $\{c_j\}_{1\le j\le k}$ on a bumpy Finsler $n$-sphere
$(S^n,F)$ with $\hat{i}(c_j)>0$. Then the following identity holds
\be \sum_{1\le j\le k,\; 1\le m\le N_j}
       (-1)^{i(c_j^m)}k_0^{\ep}(c_j^m)\frac{1}{N_j\hat{i}(c_j)}
  = \left\{\matrix{-\frac{n}{2(n-1)}, & {\rm if\;}$n$\; {\rm is\;even}, \cr
                  \frac{n+1}{2(n-1)}, & {\rm if\;}$n$\; {\rm is\; odd}. \cr}\right.,
           \lb{5.5}\ee
where $N_j=N(c_j)\in\N$ is the number defined in Lemma 5.2 for
$c_j$, $k_0^{\ep}(c_j^m)$s are the critical type numbers of
$c_j^m$, $\ep\equiv\ep(c_j^m)=(-1)^{i(c_j^m)-i(c_j)}$.  }

{\bf Proof.} Because $\dim \ol{C}_q(E,c_j^m)$ can be non-zero only
for $q=i(c_j^m)$ by Proposition 2.1, the formal Poincar\'e series
$M(t)$ becomes \be M(t)=\sum_{1\le j\le k,\; m\ge 1}
              k_0^{\ep}(c^m_j)t^{i(c_j^m)}
  = \sum_{1\le j\le k,\; 1 \le m\le N_j,\; s\ge 0}
               k_0^{\ep}(c_j^m)t^{i(c_j^{sN_j+m})}, \lb{5.6}\ee
where the last equality follows from Lemma 5.2. Write
$M(t)=\sum_{h=0}^{\infty}w_ht^h$. Then we have \be w_h\ =
\sum_{1\le j\le k,\; 1 \le m\le N_j}
              k_0^{\ep}(c_j^m)\,^\#\{s\in\N_0\,|\,i(c_j^{sN_j+m})=h\},
     \lb{5.7}\ee
where ${}^\#A$ denotes the total number of elements in a set $A$.

{\bf Claim 1.} {\it $\{w_h\}_{h\ge 0}$ is bounded}.

In fact, we have
\bea  ^\#\{s\in\N_0 &|& i(c_j^{sN_j+m})=h \}\nn\\
&=&\;^\#\{s\in\N_0 \;| \; i(c_j^{sN_j+m})=h,\;
                          |i(c_j^{sN_j+m})-(sN_j+m)\hat{i}(c_j)|\le n-1\} \nn\\
&\le &\;^\#\{s\in\N_0 \;| \;|h-(sN_j+m)\hat{i}(c_j)|\le n-1\}  \nn\\
&=&\;^\#\left\{s\in\N_0 \; \left|\;\frac{}{}\right.
      \;\frac{h-n+1-m\hat{i}(c_j)}{N_j\hat{i}(c_j)}\le s
       \le \frac{h+n-1-m\hat{i}(c_j)}{N_j\hat{i}(c_j)}\right\}  \nn\\
&\le&\; \frac{2(n-1)}{N_j\hat{i}(c_j)}+1,  \nn\eea
where the first equality follows from the fact $|i(c^m)-m\hat{i}(c)|\le n-1$
(cf. Theorem 1.4 on p. 69 of \cite{Rad1}). Hence Claim 1 holds.

Next we estimate $M^n(-1)$. By (\ref{5.7}) we have
\bea M^r(-1)
&=& \sum_{h=0}^r w_h(-1)^h   \nn\\
&=& \sum_{1\le j\le k,\; 1 \le m\le N_j}
            (-1)^{i(c_j^m)}k_0^{\ep}(c_j^m)
              \,^\#\{s\in\N_0 \,|\, i(c_j^{sN_j+m})\le r\}. \lb{5.8}\eea

{\bf Claim 2.} {\it There is a real constant $C>0$ independent of
$r$, but depending on $c_j$ for $1\le j\le k$  such that }
\be \left|M^r(-1)-\sum_{1\le j\le k,\; 1\le m\le N_j}
 (-1)^{i(c_j^m)}k_0^{\ep}(c_j^m)\frac{r}{N_j\hat{i}(c_j)}\right|  \le C.
    \lb{5.9}\ee

In fact, we have
\bea ^\#\{s\in\N_0 &|& i(c_j^{sN_j+m})\le r\}   \nn\\
&=&\;^\#\{s\in\N_0 \;| \; i(c_j^{sN_j+m})\le r,\;
             |i(c_j^{sN_j+m})-(sN_j+m)\hat{i}(c_j)|\le n-1\}  \nn\\
&\le&\;^\#\{s\in\N_0 \;| \;0\le (sN_j+m)\hat{i}(c_j)\le r+n-1\}  \nn\\
&=&\;^\#\left\{s\in\N_0 \; \left |\;\frac{}{}\right.
    \;0\le s\le \frac{r+n-1-m\hat{i}(c_j)}{N_j\hat{i}(c_j)}\right\}  \nn\\
&\le&\; \frac{r+n-1}{N_j\hat{i}(c_j)},  \lb{5.10}\eea
where the last inequality uses $\frac{1}{2}\le\frac{m}{N_j}\le 2$ by the
definition of $N_j$ and $1\le m\le N_j$.

On the other hand, we have
\bea  ^\#\{s\in\N_0 &|& i(c_j^{sN_j+m})\le r\}  \nn\\
&=&\;^\#\{s\in\N_0 \;| \; i(c_j^{sN_j+m})\le r,\;
               |i(c_j^{sN_j+m})-(sN_j+m)\hat{i}(c_j)|\le n-1\}  \nn\\
&\ge&\;^\#\{s\in\N_0\;|\;i(c_j^{sN_j+m})\le(sN_j+m)\hat{i}(c_j)+(n-1)\le r\} \nn\\
&\ge&\;^\#\left\{s\in\N_0 \; \left |\;\frac{}{}\right.
   \;0\le s\le \frac{r-n+1-m\hat{i}(c_j)}{N_j\hat{i}(c_j)}\right\}  \nn\\
&\ge&\;\frac{r-n+1}{N_j\hat{i}(c_j)}-2.  \lb{5.11}\eea
where the last inequality uses $\frac{1}{2}\le\frac{m}{N_j}\le 2$ by the
definition of $N_j$ and $1\le m\le N_j$.

By (\ref{5.10}) and (\ref{5.11}), we obtain (\ref{5.9}).

Since the sequence $\{w_h\}$ is bounded and $w_r=b_r+q_r+q_{r-1}$.
the sequence $\{q_h\}_{h\ge 0}$ of $Q(t)$ is bounded. Hence, by
(\ref{5.3}) we obtain \be \lim_{r\to\infty}\frac{1}{r}M^r(-1)
          =\lim_{r\to\infty}\frac{1}{r}P^r(-1)
  = \left\{\matrix{-\frac{n}{2(n-1)}, & {\rm if\;}$n$\; {\rm is\;even}, \cr
                  \frac{n+1}{2(n-1)}, & {\rm if\;}$n$\; {\rm is\; odd}. \cr}\right. \ee
Hence (\ref{5.5}) holds.  \hfill\hb

{\bf Remark 5.4.} V. Bangert and Y. Long in \cite{BaL1} as well as Y. Long and W. Wang
in \cite{LoW1} established such an mean index identity with exact coefficients on
Finsler $2$-spheres. For readers convenience, following ideas in \cite{BaL1} and
\cite{LoW1} we give a complete proof for $(S^n,F)$ here.

\setcounter{equation}{0}
\section{Proof of Theorem 1.2}%{Section 6}

Assuming the contrary, we prove Theorem 1.2 by contradiction. That is,
assume the following condition in this section:

{\bf (F) There exists only one prime closed geodesic $c$ on
the bumpy Finsler $S^n=(S^n,F)$.}

{\bf Lemma 6.1.} {\it Under the assumption (F), the mean index of
the closed geodesic $c$ must satisfy $\hat{i}(c)> 0$. }

{\bf Proof.} If $\hat{i}(c)=0$, then we have $i(c^m)=0$ for all $m\ge
1$ (cf. Corollary 4.2 of \cite{LLo1}). By Theorem 3.1, we have
$b_{n-1}=1$. By Proposition 2.1,we have
$$ \overline{C}_0( E,c^m)  = \Q,\quad \overline{C}_q( E,c^m) = 0
          \quad\mbox{for}\;\;q\in \N.  $$
By Theorem 3.2, we have $0=M_{n-1}\ge b_{n-1}=1$ with $n\ge 2$,
which implies $\hat{i}(c)>0$. \hfill\hb

{\bf Lemma 6.2.} {\it Under the assumption (F), the index of the
closed geodesic $c$ must satisfy $i(c)\le n-1$. }

{\bf Proof.} By contradiction, assume $i(c)>n-1$. By Corollary 4.2
of \cite{LLo1} (cf. (i) of Theorem 12.1.1 of \cite{Lon2}), we have
\be i(c^m)\ge i(c),\; \forall\; m\in \N. \lb{6.1}\ee
Hence $i(c^m)>n-1$ for all $m\ge 1$. By Proposition 2.1,we have
$$ \overline{C}_q(E,c^m) = 0 \quad\mbox{for}\;\;q\in [0,n-1]\cap\N_0. $$
Hence $\overline{C}_{n-1}(E,c^m)=0$ for all $m\in\N$, which implies
$M_{n-1}=0$. By Theorem 3.2, we have $0=M_{n-1}\ge b_{n-1}=1$, which is
a contradiction. So we have $i(c)\le n-1$. \hfill\hb

{\bf Lemma 6.3.} {\it Under the assumption (F), the index of the
closed geodesic $c$ must satisfy $$i(c)\ge n-1,$$ if one of the
following conditions is satisfied:

(i) When $n$ is even, $i(c)$ is odd and the Morse-type numbers $M_{2k}=0$
for all $k\in\N$;

(ii) When $n$ is odd, $i(c)$ is even and the Morse-type numbers
$M_{2k-1}=0$ for all $k\in\N$.}

{\bf Proof.} Suppose (i) is satisfied. Assume $i(c)<n-1$. Then
$i(c)=2k_0-1\le n-3$ for some $k_0\in\N$. So the Morse-type number
$M_{2k_0-1}\ge 1$ by Propositions 2.1 and 2.3. By (\ref{6.1}) and
$i(c)=2k_0-1$, we obtain $M_{2k-1}=0$ for $0\le k\le k_0-1$. But
by (i) of Theorem 3.1, $b_{k}=0$ for any ${k}<n-1$. Therefore, by
Theorem 3.2 and the condition (i), we have \be -1\ge-M_{2k_0-1}
=M_{2k_0} - M_{2k_0-1} + \cdots -M_{0}
  \ge b_{2k_0} - b_{2k_0-1}+\cdots -b_{0}=0, \lb{6.2}\ee
which is a contradiction.

Suppose (ii) is satisfied. Assume $i(c)<n-1$. Then $i(c)=2k_0\le
n-3$ for some $k_0\in\N$. So the Morse-type number $M_{2k_0}\ge 1$
by Propositions 2.1 and 2.3, and $b_{{k}}=0$ for any ${k}<n-1$ by
(i) of Theorem 3.1. Therefore, similarly to (\ref{6.2}), by
Theorem 3.2 we have
\be -1\ge-M_{2k_0} =M_{2k_0+1} - M_{2k_0} + \cdots -M_{0}
  \ge b_{2k_0+1} - b_{2k_0}+ \cdots -b_{0}=0, \lb{6.3}\ee
which is a contradiction. The Lemma 6.3 is proved. \hfill\hb

As an immediate consequence of Lemmas 6.2 and 6.3, we have the
following Corollary.

{\bf Corollary 6.4.} {\it Under the conditions (F), and (i) or
(ii) of Lemma 6.3, the index of the closed geodesic $c$ must
satisfy $i(c)= n-1$. }

By Theorem 3.1, the first appearance of the Betti number $b_q$ which takes
the value $2$ is when
$$ q=\left\{\matrix{
    3(n-1), \quad {\rm if}\;n\;{\rm is\;even},\cr
    2(n-1), \quad {\rm if}\;n\;{\rm is\;odd}.\cr}\right. $$
To study Case 1 in Step 1 and Case 1 in Step 2 in our proof of Theorem 1.2 below,
a basic idea is that we want to find a contradiction in both Cases by Morse
inequality before the Betti number $b_q$ reaching the value $2$ , i.e., when
$$ q\le\left\{\matrix{
    3(n-1)-2, \quad {\rm if}\;n\;{\rm is\;even},\cr
    2(n-1)-2, \quad {\rm if}\;n\;{\rm is\;odd}.\cr}\right. $$
Hence we construct the following sets
\be \Theta(n)=\left\{\matrix{
    \{j\in 2\N-1\ |\ n-1\le j\le 3n-5\}, \quad {\rm when}\;n\;{\rm is\;even},\cr
    \{j\in 2\N\ |\ n-1\le j\le 2n-4\}, \quad {\rm when}\;n\;{\rm is\;odd}.\cr}\right.
    \lb{6.4}\ee

{\bf Lemma 6.5.} {\it Assume the conditions (F), and (i) or (ii) of Lemma
6.3 hold. Suppose $\frac{i(c^m)-i(c)}{2}\in\{0,1,\cdots,m-1\}$ for all $m$.
Then for every $n-1+2k\in \Theta(n)$, there exists a unique iteration
$c^{k+1}$ of $c$ such that }
\be i(c^{k+1})=n-1+2k. \lb{6.5}\ee

{\bf Proof.} By Corollary 6.4, $i(c)=n-1$. So (\ref{6.5}) holds
when $k=0$. Assume that there exists another iteration $c^{m_0}$
such that $i(c^{m_0})=i(c)=n-1$, by Proposition 2.1 and
$i(c^m)-i(c)\in 2\N_0$ for all $m$, it yields $M_{n-1}\ge 2$ and
$M_r=0$ for any $r<n-1$ by (\ref{6.1}). On the other hand, by
Theorem 3.1, $b_{n-1}=1$ and $b_n=b_r=0$ for any $r<n-1$. Noting
that $M_n=0$ by the condition (i) or (ii), by Theorem 3.2 we have
$$ -2 \ge M_{n} - M_{n-1} + \cdots +(-1)^n M_{0}
      \ge b_{n} - b_{n-1}+ \cdots +(-1)^n b_{0}=-1, $$
which is a contradiction. Thus Lemma 6.5 holds for $k=0$.

By induction, we assume that Lemma 6.5 holds for all $k\le\hat{k}$,
where $1\le\hat{k}<\max\Theta(n)$, i.e., there exists a unique iteration
$c^{k+1}$ such that (\ref{6.5}) holds for all $k\le\hat{k}$. Hence
we have
\be M_q=1\quad\mbox{for}\quad
  q\in\hat{\Theta}(n)\equiv\{n-1+2t\ |\ t\in
     [0,\hat{k}]\cap\N_0\}\subset \Theta(n),  \lb{6.6}\ee
by Proposition 2.1 and $i(c^m)-i(c)\in 2\N_0$ for all $m$.

Firstly, we prove that $i(c^{\hat{k}+2})=n-1+2(\hat{k}+1)$. By
the condition $\frac{i(c^m)-i(c)}{2}\in\{0,1,\cdots,m-1\}$ for all
$m$, it yields $i(c^{\hat{k}+2})=n-1+2t,t\in\{0,1,\cdots,\hat{k}+1\}$.
If $t\in\{0,1,\cdots,\hat{k}\}$, then there must exist some
$s\in[0,\hat{k}]\cap\N_0$ such that $i(c^{\hat{k}+2})=i(c^{s+1})=n-1+2s$,
which contradicts to the uniqueness of $c^{s+1}$. So the only possibility
is $i(c^{\hat{k}+2})=n-1+2(\hat{k}+1)$.

On the other hand, assume that there exists another iteration $c^{m_0}$
such that $i(c^{m_0})=i(c^{\hat{k}+2})=n-1+2(\hat{k}+1)\equiv\ka\in \Theta(n)$.
Then it yields $M_{\ka}\ge 2$ by Proposition 2.1 and $i(c^m)-i(c)\in 2\N_0$
for all $m$. Note that for $l\in[0,\ka-1]\setminus\hat{\Theta}(n)$, there
holds $M_{\ka+1}=M_l=0$ by the condition (i) or (ii). Therefore among all the
$M_l$'s with $l<\ka$, half of them are zero, and half of them are $1$. By
Theorem 3.1, $b_{\ka+1}=b_l=0$ for the same $l$ mentioned above
and $b_{\ka}=b_q=1$ for $q$ in (\ref{6.6}). Note that $\ka-n\in (2\Z+1)$ by
the definition of $\ka$. Therefore, by Theorem 3.2 we have
\bea   -(2+\frac{\ka-n+1}{2})
&\ge& M_{\ka+1} - M_{\ka}+M_{\ka-1}+ \cdots -M_{n-1} \nn\\
&\ge& b_{\ka+1} - b_{\ka}-b_{\ka-1}+ \cdots -b_{n-1} \nn\\
&=& -(1+\frac{\ka-n+1}{2}), \nn\eea
which is a contradiction. Therefore Lemma 6.5 holds for $k=\hat{k}+1$. This
completes the proof. \hfill\hb

Now we can give

{\bf Proof of Theorem 1.2.} We carry out the proof in two steps
under the assumption (F) on $(S^n,F)$.

{\bf Step 1.} {\it When $n$ is even, we claim that there must be another prime
closed geodesic.}

We study the problem in five cases according to our classification in
Section 4.

{\it Case 1.} {\it $c$ belongs to NCG-1.}

In this case,  $i(c^m)$ is odd and $i(c^m)-i(c)$ is even for any
$m\in\N$. Hence $\ep=\ep(c^m)=1$ and $k_0^{\ep}(c^m)=1$ by
Proposition 2.3. By Lemma 5.2 and Theorem 5.3, we have
$-\frac{1}{\hat{i}(c)}=-\frac{n}{2(n-1)}$, i.e.,
\be 2p+\sum_{i=1}^{n-2r-1}\frac{\th_i}{\pi}=\hat{i}(c)=
               \frac{2(n-1)}{n}<2.  \lb{6.7}\ee
On the other hand, by Proposition 2.1 and the oddness of $i(c^m)$,
we have the Morse-type numbers \be M_k=0\qquad \mbox{for all}\;
k\in 2\N.   \lb{6.8}\ee Hence, by Corollary 6.4, $2p+n-2r-1=i(c)=
n-1$, which yields $p=r\in\N_0$. Together with (\ref{6.7}), it
yields $r=p=0$. So in this case, we have \be i(c)=n-1,\quad
i(c^m)=2\sum_{i=1}^{n-1}\left[\frac{m\th_i}{2\pi}\right]+n-1\quad
\mbox{and}\quad \hat{i}(c)=\sum_{i=1}^{n-1}\frac{\th_i}{\pi} =
\frac{2(n-1)}{n}<2. \lb{6.9}\ee Note that, by (\ref{6.9}), for any
$m\in\N$, there holds \be
\sum_{i=1}^{n-1}\frac{m\th_i}{2\pi}=\frac{m(n-1)}{n}<m.\lb{6.10}\ee
Hence \be  \sum_{i=1}^{n-1}\left[\frac{m\th_i}{2\pi}\right]\in
\{0,1,\cdots,m-1\}.
          \lb{6.11}\ee

{\it Claim 1: For any $m\in [1, n-1]\cap\N$, there holds $i(c^m)=2(m-1)+n-1$.}

In fact, we have $i(c)=n-1$ by (\ref{6.9}). By induction, assume
that $i(c^k)=2(k-1)+n-1$ for some $k\in[1,n-2]\cap\N$. By
(\ref{6.9}) and (\ref{6.11}), for $k+1$ we have
\be i(c^{k+1})=2\sum_{i=1}^{n-1}\left[\frac{(k+1)\th_i}{2\pi}\right]+n-1\in
      \{n-1+2s\,|\,s\in [0,k]\cap\N_0\}.  \lb{6.12}\ee
Note that the conditions of Lemma 6.5 are satisfied by (\ref{6.8}), (\ref{6.9})
and (\ref{6.11}). If $i(c^{k+1})\le n-1+2(k-1)$, then $i(c^{k+1})=i(c^r)$ for
some $r\in [1,k]\cap\N$, which yields a contradiction to Lemma 6.5. So the
only possibility is $i(c^{k+1})=2k+n-1$. Claim 1 is proved.

Next we consider $i(c^{n})$. By (\ref{6.10}), we have
\be \sum_{i=1}^{n-1}\frac{n\th_i}{2\pi}=\frac{n(n-1)}{n}=n-1. \lb{6.13}\ee
Noting that
$\{\frac{n\th_i}{2\pi}\}$ is an irrational number sequence, by
(\ref{6.13}) it yields
\be \sum_{i=1}^{n-1}\left[\frac{n\th_i}{2\pi}\right]
            \in\{0,1,\cdots,n-2\}.\lb{6.14}\ee
Hence
\be i(c^n)=2\sum_{i=1}^{n-1}\left[\frac{n\th_i}{2\pi}\right]+n-1\in
    \{n-1+2s\,|\, s\in [0,n-2]\cap\N_0\}. \lb{6.15}\ee
Therefore, by Claim 1 and (\ref{6.15}), there must exists some
$r\in[1,n-1]\cap\N$ such that $i(c^n)=i(c^r)$, which is also a
contradiction to Lemma 6.5.

{\it Case 2. $c$ belongs to NCG-2.}

{\it Subcase 2.1.} If $i(c)=p$ is even, we have $N=1$ by Lemma 5.2
and $i(c^m)$ is even for any $m\in\N$. Hence $\ep=1$ and
$k_0^{\ep}(c^m)=1$ by Proposition 2.3. Thus, by Theorem 5.3, it
yields $\frac{1}{\hat{i}(c)}=-\frac{n}{2(n-1)}$, which yields a
contradiction by Lemma 6.1.

{\it Subcase 2.2.} If $i(c)=p$ is odd, $N=2$ by Lemma 5.2 and
\be  i(c^m)\; \mbox{is}\;
\left\{\matrix{\mbox{even}, & {\rm if\;}$m$\; {\rm is\; even}, \cr
               \mbox{odd}, & {\rm if\;}$m$\; {\rm is\; odd}. \cr}\right.
\lb{6.16}\ee
So by Theorem 5.3, it yields $\frac{1}{\hat{i}(c)}=\frac{n}{n-1}$,
which implies
\be \hat{i}(c)=p-k+\sum_{i=1}^{k}\frac{\th_i}{\pi}=\frac{n-1}{n}<1.
    \lb{6.17}\ee
On the other hand, by (\ref{6.16}) and Proposition 2.1, we have
the Morse-type numbers $M_k=0$ for any $k\in 2\N$. So by Corollary
6.4, we have $i(c)=p=n-1$. Hence by (\ref{6.17}), it yields $n-2<k$.
This yields a contradiction to $k\le n-2r-2$ in the definition
of NCG-2.

{\it Case 3. $c$ belongs to NCG-3.}

{\it Subcase 3.1.} If $i(c)=p$ is even, we obtain $N=2$ by Lemma 5.2
and $i(c^2)$ is odd. Hence $\ep=\ep(c^2)=-1$ and $k_0^{\ep}(c^2)=0$
So by Theorem 5.3, it yields $\frac{1}{2\hat{i}(c)}=-\frac{n}{2(n-1)}$,
which yields a contradiction by Lemma 6.1.

{\it Subcase 3.2.} If $i(c)=p$ is odd, we obtain $N=1$ by Lemma 5.2 and
$i(c^m)$ is odd for any $m$. Hence $\ep=1$ and $k_0^{\ep}(c^m)=1$
by Proposition 2.3. So by Theorem 5.3, it yields
$-\frac{1}{\hat{i}(c)}=-\frac{n}{2(n-1)}$, which implies
\be \hat{i}(c)=p-k+\sum_{i=1}^{k}\frac{\th_i}{\pi}=\frac{2(n-1)}{n}<2.
   \lb{6.18}\ee
Noting that $p-k$ is even, this yields $p\le k$. Because the
Morse-type numbers $M_k=0$ holds for any $k\in 2\N$ by Proposition
2.1, by Corollary 6.4 we have $i(c)=p=n-1$. Hence $n-1 = p\le k$.
This yields a contradiction to $k\le n-2r-2$ in the definition of
NCG-3.

{\it Case 4. $c$ belongs to NCG-4.}

{\it Subcase 4.1.} If $i(c)=p$ is even, $N=2$ by Lemma 5.2 and $i(c^2)$
is odd. Hence $\ep=\ep(c^2)=-1$ and $k_0^{\ep}(c^2)=0$. So by Theorem 5.3,
it yields $\frac{1}{2\hat{i}(c)}=-\frac{n}{2(n-1)}$, which yields a
contradiction by Lemma 6.1.

{\it Subcase 4.2.} If $i(c)=p$ is odd, $N=1$ by Lemma 5.2 and
$i(c^m)$ is odd for all $m$. Hence $\ep=1$ and $k_0^{\ep}(c^m)=1$
by Proposition 2.3. Thus by Theorem 5.3, it yields
$\frac{1}{\hat{i}(c)}=\frac{n}{2(n-1)}$. But in this case
$\hat{i}(c)=(p-1)+\frac{\th_1}{\pi}$ is an irrational number. This
leads to a contradiction.

{\it Case 5. $c$ belongs to NCG-5.}

{\it Subcase 5.1.} If $i(c)=p$ is even, $N=1$ by Lemma 5.2 and
$i(c^m)$ is even for all $m$. Hence $\ep=1$ and $k_0^{\ep}(c^m)=1$
by Proposition 2.3. So by Theorem 5.3, there holds
$\frac{1}{\hat{i}(c)}=-\frac{n}{2(n-1)}$, which yields a
contradiction by Lemma 6.1.

{\it Subcase 5.2.} If $i(c)=p$ is odd, we obtain $N=2$ by Lemma
5.2 and $i(c^2)$ is even. Hence $\ep=\ep(c^2)=-1$ and
$k_0^{\ep}(c^2)=0$. Thus by Theorem 5.3, it yields
$\frac{1}{\hat{i}(c)}=\frac{n}{n-1}$. Hence
$\hat{i}(c)=p=\frac{n-1}{n}$, which yields a contradiction because
$p\in\N_0$.

Therefore when $n$ is even, we have proved that there must exist another
prime closed geodesic on the bumpy Finsler $n$-sphere $(S^n,F)$.

{\bf Step 2.} {\it When $n$ is odd, we claim that there must be another
prime closed geodesic.}

We continue our study in five cases according to the classification in
Section 4.

{\it Case 1. $c$ belongs to NCG-1.}

In this case, $i(c^m)$ is even for all $m$. Hence $\ep=1$ and
$k_0^{\ep}(c^m)=1$ by Proposition 2.3. By Lemma 5.2 and Theorem
5.3, we have $\frac{1}{\hat{i}(c)}=\frac{n+1}{2(n-1)}$, i.e.,
\be 2p+\sum_{i=1}^{n-2r-1}\frac{\th_i}{\pi}=\hat{i}(c)
          =\frac{2(n-1)}{n+1}<2.   \lb{6.19}\ee
On the other hand, by Proposition 2.1 and the evenness of
$i(c^m)$, we have the Morse-type numbers \be M_k=0 \qquad
\forall\, k\in 2\N-1.   \lb{6.20}\ee Thus by Corollary 6.4, we
have $2p+n-2r-1=i(c)=n-1$, which implies $p=r\in\N_0$. Together
with (\ref{6.19}), it yields $r=p=0$. Hence in this case, we have
\be i(c)=n-1,\quad
i(c^m)=2\sum_{i=1}^{n-1}\left[\frac{m\th_i}{2\pi}\right]+n-1\quad
\mbox{and}\quad\hat{i}(c)=\sum_{i=1}^{n-1}\frac{\th_i}{\pi}=
\frac{2(n-1)}{n+1}<2.      \lb{6.21}\ee

Note that, by (\ref{6.21}), for any $m\in\N$ there holds
\be \sum_{i=1}^{n-1}\frac{m\th_i}{2\pi}=\frac{m(n-1)}{n+1}<m. \lb{6.22}\ee
Hence
\be \sum_{i=1}^{n-1}\left[\frac{m\th_i}{2\pi}\right]\in \{0,1,\cdots,m-1\}.
    \lb{6.23}\ee
Let
\be m_1=\frac{n-1}{2},\; m_2=\frac{n+1}{2}.\lb{6.24}\ee

{\it Claim 1: For any $m\in [1, m_1]\cap\N$, there holds
$i(c^m)=2(m-1)+n-1$.}

In fact, we have $i(c)=n-1$. By induction, assume that
$i(c^k)=2(k-1)+n-1$ for some $k\in[1,m_1-1]\cap\N$. By
(\ref{6.21}) and (\ref{6.23}), for $k+1$ we have
\be  i(c^{k+1})=2\sum_{i=1}^{n-1}\left[\frac{(k+1)\th_i}{2\pi}\right]+n-1
       \in \{n-1+2s\,|\, s\in [0,k]\cap\N_0\}.    \lb{6.25}\ee
Note that the conditions of Lemma 6.5 are satisfied by (\ref{6.20}),
(\ref{6.21}) and (\ref{6.23}). If $i(c^{k+1})\le n-1+2(k-1)$, then
$i(c^{k+1})=i(c^r)$ holds for some $r\in [1,k]\cap\N$, which yields a
contradiction to Lemma 6.5. Therefore the only possibility is
$i(c^{k+1})=2k+n-1$, which proves Claim 1.

Next we consider $i(c^{m_2})$. By (\ref{6.22}) and (\ref{6.24}),
we have
\be \sum_{i=1}^{n-1}\frac{m_2\th_i}{2\pi}=
        \frac{(n+1)(n-1)}{2(n+1)}=\frac{n-1}{2}.  \lb{6.26}\ee
Noting that $\{\frac{m_2\th_i}{2\pi}\}$ is an irrational number
sequence, by (\ref{6.26}) it yields \be
\sum_{i=1}^{n-1}\left[\frac{m_2\th_i}{2\pi}\right]
    \in\{0,1,\cdots,\frac{n-1}{2}-1\}.  \lb{6.27}\ee
Hence
\be i(c^{m_2})=2\sum_{i=1}^{n-1}\left[\frac{m_2\th_i}{2\pi}\right]+n-1
    \in \{n-1+2s \,|\, s\in [0, \frac{n-1}{2}-1]\cap\N_0\}. \lb{6.28}\ee
Therefore, by Claim 1 and (\ref{6.28}), there must exist some
$r\in[1,m_1]\cap\N$ such that $i(c^{m_2})=i(c^r)$, which is also a
contradiction to Lemma 6.5.

{\it Case 2. $c$ belongs to NCG-2.}

{\it Subcase 2.1.} If $i(c)=p$ is even, $N=1$ by Lemma 5.2 and
$i(c^m)$ is even for all $m$. Hence $\ep=1$ and $k_0^{\ep}(c^m)=1$
by Proposition 2.3. So by Theorem 5.3, it yields
$\frac{1}{\hat{i}(c)}=\frac{n+1}{2(n-1)}$, which implies
\be \hat{i}(c)=p-k+\sum_{i=1}^{k}\frac{\th_i}{\pi}
     =\frac{2(n-1)}{n+1}<2. \lb{6.29}\ee
Noting that $p-k$ is even, this yields $p\le k$. Because the
Morse-type numbers $M_k=0$ for all $k\in 2\N-1$, by Corollary 6.4
we have $i(c)=p= n-1$. Hence $n-1 =p\le k$. This contradicts to
$k\le n-2r-2$ by the definition of NCG-2.

{\it Subcase 2.2.} If $i(c)=p$ is odd, we have $N=2$ by Lemma 5.2
and $i(c^2)$ is even. Hence $\ep=\ep(c^2)=-1$ and
$k_0^{\ep}(c^2)=0$. So by Theorem 5.3, it yields
$-\frac{1}{\hat{i}(c)}=\frac{n+1}{n-1}$, which yields a
contradiction.

{\it Case 3. $c$ belongs to NCG-3.}

{\it Subcase 3.1.} If $i(c)=p$ is even, $N=2$ by Lemma 5.2 and
\be i(c^m)\; \mbox{is}\;
\left\{\matrix{\mbox{even}, & {\rm if\;}$m$\; {\rm is\; odd}, \cr
               \mbox{odd}, & {\rm if\;}$m$\; {\rm is\; even}. \cr}\right.
    \lb{6.30}\ee
Therefore by Theorem 5.3, it yields $\frac{1}{\hat{i}(c)}=\frac{n+1}{n-1}$,
which implies
\be \hat{i}(c)=p-k+\sum_{i=1}^{k}\frac{\th_i}{\pi}=\frac{n-1}{n+1}<1.
   \lb{6.31}\ee
By (\ref{6.30}) and Proposition 2.1, we have the Morse-type
numbers $M_k=0$ for all $k\in 2\N-1$. Thus by Corollary 6.4, we
have $i(c)=p=n-1$. Hence, by (\ref{6.31}), it yields $n-2< k$.
This is a contradiction to $k\le n-2r-2$ by the definition of
NCG-3.

{\it Subcase 3.2.} If $i(c)=p$ is odd, $N=1$ by Lemma 5.2. So, by
Theorem 5.3, it yields $-\frac{1}{\hat{i}(c)}=\frac{n+1}{2(n-1)}$,
which yields a contradiction by Lemma 6.1.

{\it Case 4. $c$ belongs to NCG-4.}

{\it Subcase 4.1.} If $i(c)=p$ is even, $N=2$ by Lemma 5.2 and
$i(c^2)$ is odd. Hence $\ep=\ep(c^2)=-1$ and $k_0^{\ep}(c^2)=0$.
So by Theorem 5.3, it yields $\frac{1}{\hat{i}(c)}=\frac{n+1}{n-1}$.
But in this case $\hat{i}(c)=(p-1)+\frac{\th_1}{\pi}$ is an
irrational number. This leads to a contradiction.

{\it Subcase 4.2.} If $i(c)=p$ is odd, we have $N=1$ by Lemma 5.2.
Thus by Theorem 5.3, we have $-\frac{1}{\hat{i}(c)}=\frac{n+1}{2(n-1)}$,
which yields a contradiction by Lemma 6.1.

{\it Case 5. $c$ belongs to NCG-5.}

{\it Subcase 5.1.} If $i(c)=p$ is even, $N=1$ by Lemma 5.2. So by
Theorem 5.3, it yields
$\frac{1}{p}=\frac{1}{\hat{i}(c)}=\frac{n+1}{2(n-1)}$. Hence we
have \be1>\frac{n-1}{n+1}=\frac{p}{2}\ge 1,\ee which is a
contradiction.

{\it Subcase 5.2.} If $i(c)=p$ is odd, $N=2$ by Lemma 5.2 and
$i(c^2)$ is even. Hence $\ep=\ep(c^2)=-1$ and $k_0^{\ep}(c^2)=0$.
Thus by Theorem 5.3, there holds
$-\frac{1}{\hat{i}(c)}=\frac{n+1}{n-1}$, which yields a
contradiction by Lemma 6.1.

So when $n$ is odd, there must exist another prime closed geodesic
on the bumpy Finsler $n$-sphere $(S^n,F)$.

The Steps 1 and 2 complete the proof of Theorem 1.2. \hfill\hb

\medskip

\noindent {\bf Acknowledgements.} The authors sincerely thank the
referee for his/her valuable comments on this paper.

\bibliographystyle{abbrv}

%\bigskip

%\noindent  Huagui Duan,

%Chern Institute of Mathematics, Nankai University, Tianjin 300071,
%the People's Republic of China.
%Email: duanhuagui@163.com

%\bigskip

%\noindent Yiming Long,

%hern Institute of Mathematics and LPMC, Nankai University, Tianjin 300071,
%the People's Republic of China.
%Email: longym@nankai.edu.cn

\end{document}